\documentclass[preprint,12pt]{elsarticle}

\usepackage{amssymb}
\usepackage{amscd,amsmath,amsfonts,graphicx}
\usepackage{amsthm}

\begin{document}
\baselineskip=22pt \centerline{\large \bf Characterizations of the
Normal Distribution via the Independence} \centerline{\large \bf  of
the Sample Mean and the Feasible Definite Statistics}
 \centerline{\large \bf  with
Ordered Arguments}

\vspace{1cm} \centerline{Chin-Yuan Hu $^a$ and Gwo Dong Lin
$^{b,c}$} \centerline{$^a$ National Changhua University of
Education, Taiwan} \centerline{$^b$ Hwa-Kang Xing-Ye Foundation,\
Taiwan} \centerline{$^c$ Academia Sinica, Taiwan}

 \vspace{1cm} \noindent {\bf Abstract.} It is well known that the independence
 of the sample mean and the sample variance characterizes the normal distribution.
By using Anosov's theorem, we further investigate the analogous
characteristic properties in terms of the sample mean and some
feasible definite statistics. The latter statistics introduced in
this paper for the first time are based on nonnegative, definite and
continuous functions of ordered arguments with positive degree of
homogeneity. The proposed approach seems to be natural and can be
used to derive easily characterization results for many feasible
definite statistics, such as known characterizations involving the
sample variance, sample range as well as Gini's mean difference.
\vspace{0.1cm}\\
\hrule
{\bf Running title: Characterizations of the
Normal Distribution }\\
{\bf Corresponding Author:} Gwo Dong Lin (gdlin@stat.sinica.edu.tw)\\
{\bf Postal addresses:} Chin-Yuan Hu, National Changhua University
of Education, No.\,1, Jinde Road, Changhua 50007, Taiwan (ROC). (E-mail: buhuua@gmail.com)\\
Gwo Dong Lin, (1) Social and Data Science Research Center, Hwa-Kang
Xing-Ye Foundation, No.\,16, Lane 276, Rui'an Street, Da'an
District, Taipei 10659, Taiwan (ROC) and (2) Institute of
Statistical Science, Academia Sinica, No.\,128, Section 2, Academia
Road, Nankang District, Taipei 11529, Taiwan (ROC). (E-mail:
 gdlin@stat.sinica.edu.tw)
\newpage
\noindent{\bf 1. Introduction}
\newcommand{\bin}[2]{
   \left(
     \begin{array}{@{}c@{}}
         #1  \\  #2
     \end{array}
   \right)          }

In the characterization theory of probability or statistical
distributions, one of the remarkable results is the characterization
of the normal distribution through the independence of the sample
mean and the sample variance described below.

Let $X$ be a random variable with distribution $F$ on the whole real
line $\mathbb R:=(-\infty,\infty),$  denoted $X\sim F.$ Let $X_1,
X_2,\ldots, X_n$ be a random sample of size $n\ge 2$ from
distribution $F.$ Denote the sample mean
$\overline{X}_n=\frac1n\sum_{i=1}^nX_i$ and the sample variance
$S_n^2=\frac{1}{n-1}\sum_{i=1}^n(X_i-\overline{X}_n)^2.$ If $F$ is
normal, then  $\overline{X}_n$ and $(X_1-\overline{X}_n,
X_2-\overline{X}_n,\ldots, X_n-\overline{X}_n)$ are independent
(see, e.g., Rohatgi 1976, p.\,321), and hence so are
$\overline{X}_n$ and $S_n^2.$ Conversely, if $\overline{X}_n$ and
$S_n^2$ are independent, then $F$ is a normal distribution
(including the degenerate case). Geary (1936) proved this result
under the extra moment condition that the underlying distribution
$F$ has finite moments of all orders and then Lukacs (1942) improved
this result to just include the second moment assumption.
 Kawata and Sakamoto (1949) as well as Zinger
(1951) solved the problem completely by using different approaches.

Instead of the sample variance $S_n^2,$ Laha (1956) considered the
quadratic form $Q=\sum_{i,j}a_{ij}X_iX_j$ and investigated the
conditions under which the independence of $\overline{X}_n$ and $Q$
characterizes the normal distribution. In this regard, see Kagan et
al.\,(1973), Chapter 4.

In this paper, we introduce the novel statistics $Z_n$ with ordered
arguments of the form:
\begin{eqnarray}Z_n=U(X_{(1)}-\overline{X}_n,X_{(2)}-\overline{X}_n,\ldots,X_{(n)}-\overline{X}_n),\end{eqnarray}
where $X_{(1)}\le X_{(2)}\le \cdots\le X_{(n)}$ are order statistics
of the random sample $X_1, X_2,\ldots, X_n$ of size $n\ge 3$  (our
approach needs to assume $n\ge 3$; see (10) and (15) below). We say
that the base function $U$ defined on the ordered set
\begin{eqnarray}{\bf A}=\Big\{(\lambda_1,\ldots,\lambda_n): \lambda_1\le \lambda_2\le\cdots\le \lambda_n, \
\small\sum_{i=1}^n\lambda_i=0\Big\}
\end{eqnarray} is  feasible definite with
positive degree $p$ of homogeneity if it satisfies
 the following conditions:\\
 (i) $U$ is nonnegative and continuous on ${\bf A},$\\
(ii) $U(\lambda_1,\ldots,\lambda_n)=0$ if and only if
$(\lambda_1,\ldots,\lambda_n)=(0,\ldots,0)$ (definiteness) and\\
(iii)
$U(s(\lambda_1,\ldots,\lambda_n))=s^pU(\lambda_1,\ldots,\lambda_n)$
for all
$s>0$ and $(\lambda_1,\ldots,\lambda_n)\in {\bf A}$ (positive degree $p$ of homogeneity).\\
The corresponding statistic $Z_n$ defined in (1) through such a base
function $U$ is called a feasible definite statistic on ${\bf A}$
with positive degree $p$ of homogeneity.

Note that the sample mean
 is in general
not a feasible definite statistic, because it may take negative
values and doesn't satisfy the definiteness condition (see (i) and
(ii) above).  We will find conditions on the base function $U$ and
the underlying distribution $F$ under which the independence of the
sample mean $\overline{X}_n$ and the feasible definite statistic
$Z_n$ characterizes the normal distribution (see the Theorem below).
On the other hand,  Hwang and Hu (2000) considered the statistics of
the form (without the concept of homogeneity):
\[Z_n=S_n\cdot \exp(\psi(\Lambda_1,\Lambda_2,\ldots,\Lambda_n)),
\]
where $\Lambda_i=(X_{(i)}-\overline{X}_n)/S_n$  and $\psi$ is a
bounded continuous real-valued function.

Compared with the above product-form statistics, our approach here
seems to be more natural and can be used to derive easily many
characterizations involving feasible definite statistics. The latter
statistics encompass the familiar sample variance $S_n^2$, sample
range
\begin{eqnarray}R_{n}:=X_{(n)}-X_{(1)}=(X_{(n)}-\overline{X}_n)-(X_{(1)}-\overline{X}_n),\end{eqnarray}
and Gini's mean difference $G_n$ (see Corollaries 1, 3 and 5 below).

 The main results are stated in Section 2. Section 3 provides the crucial tools --  Anosov's theorem, Benedetti's inequality as well as
  three lemmas.
 The proofs of the main results are given in Section 4.  Section 5 provides some more results about the sample range and
 Gini's mean difference as well as two conjectures. Finally, we have in Section 6 some discussions.
\medskip\\
\noindent{\bf 2. Main results}

Throughout the section, we consider a random variable $X\sim F$
having {\it positive continuous density} $f_X$ on $\mathbb R.$ Let
$X_{(1)}\le X_{(2)}\le \cdots\le X_{(n)}$ be the order statistics of
a random sample $X_1, X_2,\ldots, X_n$ of size $n\ge 3$ from
distribution $F.$ Then we have the following result.
\medskip\\
\noindent {\bf Theorem.} {\it Let
$Z_n=U(X_{(1)}-\overline{X}_n,X_{(2)}-\overline{X}_n,\ldots,X_{(n)}-\overline{X}_n)$
 be a feasible definite
statistic on ${\bf A}$ $($defined in $(2))$ with positive degree $p$
of homogeneity. Then $\overline{X}_n$ and $Z_n$ are independent if
and only if $F$ is normal.}

Applying the theorem to some suitable base functions yields the
following corollaries. The known results (Corollaries 1 and 5) are
listed here for comparison. In fact, they are special cases of
Corollaries 2 and 4, and hence their proofs are omitted.
\medskip\\
\noindent {\bf Corollary 1} (Hwang and Hu 2000).  {\it The sample
mean $\overline{X}_n$ and the sample range $R_{n}$ $($defined in
$(3))$ are independent if and only if $F$ is normal.}
\medskip\\
\noindent {\bf Corollary 2.} {\it  Assume that  $a_1\le a_2\le
\cdots\le a_n$ are not all equal. Then the sample mean
$\overline{X}_n$ and the feasible definite statistic
\begin{eqnarray}
Z_n=\sum_{i=1}^na_i(X_{(i)}-\overline{X}_n)\end{eqnarray} are
independent if and only if $F$ is normal.}

When $\sum_{i=1}^na_i=0,$ Corollary 2 reduces  to Corollary 2.2 of
Hwang and Hu (2000), from which the sample range $R_n$ and Gini's
mean difference $G_n$ are derived as special cases.
\medskip\\
\noindent {\bf Corollary 3.} {\it Let $p>0,$ $a_1>0, a_n>0,$ and
$a_{i}\ge 0,$ where $2\le i\le n-1.$  Then the sample mean
$\overline{X}_n$ and the feasible definite statistic\begin{eqnarray}
Z_n=\sum_{i=1}^na_{i}|X_{(i)}-\overline{X}_n|^p\end{eqnarray} are
independent if and only if $F$ is normal.}

When $p=2$ and $a_i=1/(n-1)$ for all $i,$ the $Z_n$ in (5) reduces
to the sample variance $S_n^2.$ Hence this result includes the
classical one as a special case.
\medskip\\
\noindent {\bf Corollary 4.} {\it Let $p>0,$ $a_{ij}\ge 0,$ where
$1\le i,\ j\le n,$ and $a_{1n}+a_{n1}>0.$ Then the sample mean
$\overline{X}_n$ and the feasible definite statistic\begin{eqnarray}
Z_n=\sum_{i=1}^n\sum_{j=1}^na_{ij}|X_{(i)}-X_{(j)}|^p\end{eqnarray}
are independent if and only if $F$ is normal.}
\medskip\\
\noindent {\bf Corollary 5} (Hwang and Hu 2000). {\it The sample
mean $\overline{X}_n$ and Gini's mean difference
\begin{eqnarray}G_n&=&\frac{1}{n(n-1)}\sum_{i=1}^n\sum_{j=1}^n|X_i-X_j|\nonumber\\
&=&\frac{1}{n(n-1)}\sum_{i=1}^n\sum_{j=1}^n|X_{(i)}-X_{(j)}|
=\frac{4}{n(n-1)}\sum_{i=1}^n(i-(n+1)/2)X_{(i)}
\end{eqnarray} are independent if and only if $F$ is normal.}

For the last equality in (7), see David and Nagaraja (2003),
pp.\,249 and 279.
\medskip\\
\noindent {\bf Corollary 6.} {\it Let $(a_{ij})_{i,j=1}^n$ be a
positive definite $($real\,$)$ matrix. Then the sample mean
$\overline{X}_n$ and the feasible definite statistic\begin{eqnarray}
Z_n=\sum_{i=1}^n\sum_{j=1}^na_{ij}(X_{(i)}-\overline{X}_n)(X_{(j)}-\overline{X}_n)\end{eqnarray}
are independent if and only if $F$ is normal.}
\medskip\\
\noindent {\bf Corollary 7.} {\it Let $p>0,\,q>0.$  Assume further
that $a_{ij}\ge 0,$ where $1\le i,\,j\le n,$ and $a_{11}>0,\
a_{nn}>0.$ Then the sample mean $\overline{X}_n$ and the feasible
definite statistic\begin{eqnarray}
Z_n=\sum_{i=1}^n\sum_{j=1}^na_{ij}|X_{(i)}-\overline{X}_n|^p\,|X_{(j)}-\overline{X}_n|^q\end{eqnarray}
are independent if and only if $F$ is normal.}
\medskip\\
\noindent{\bf 3. Crucial tools -- Anosov's theorem, Benedetti's
inequality and three lemmas}

To prove the main results, we need Anosov's theorem (see Anosov 1964
or Kagan et al.\,1973, Chapter 4), Benedetti's inequality (see
Benedetti 1957, 1995, Georgescu-Roegen 1959, Sarria and Martinez
2016, or Hwang and Hu 1994a), and three crucial lemmas.
\medskip\\
\noindent {\bf Anosov's Theorem.} {\it Let $n\ge 3$ be an integer
and let $X\sim F$ have a positive continuous density $f_X$ on
$\mathbb R.$ Define the $(n-2)$-dimensional torus
\begin{eqnarray}\Phi=\{\phi=(\phi_1,\ldots,\phi_{n-2}):
\phi_j\in[0,\pi],\ j=1,2,\ldots,n-3;\
\phi_{n-2}\in[0,2\pi]\}.\end{eqnarray} Let $T$ be a nonnegative
continuous function on $\Phi$ such that $\int_{\Phi}T(\phi)\,{\rm
d}\phi\in(0,\infty).$ Assume further that $\sigma_j,\ j=1,\ldots,
n,$ are continuous functions on $\Phi$ satisfying
\begin{eqnarray} \sum_{j=1}^n\sigma_j(\phi)=0\ \ {\rm  and }\
\sum_{j=1}^n\sigma_j^2(\phi)\ \in(0,\infty)\ \ {\rm for\ all }\ \
\phi\in\Phi.\end{eqnarray} Under these conditions, if
 the density $f_X$ satisfies the integro-functional equation\,$:$
\begin{eqnarray}
\int_{\Phi}\prod_{j=1}^nf_X(t+s\sigma_j(\phi))\,T(\phi)\,{\rm d}\phi=c\,(f_X(t))^n\!\!
\int_{\Phi}\prod_{j=1}^nf_X(s\sigma_j(\phi))\,T(\phi)\,{\rm d}\phi\ \ {\forall}\ t\in{\mathbb R}\ {\rm  and}\ s\ge 0,~
\end{eqnarray}
where $c>0$ is a constant, then $f_X$ is normal.}

\noindent{\bf Remark 1.} The proof of Anosov's Theorem is
complicated. However, as noted by Kagan et al.\,(1973), p.\,143, if
we further assume the density function $f_X$ in (12) to be  {\it
continuously twice-differentiable}, then the proof becomes much
simpler,  just by differentiating (12) twice at $s=0$ and solving
the obtained functional equation.
\medskip\\
\noindent{\bf Benedetti's Inequality.} {\it Let $n\ge 2$ be an
integer and let ${\bf \mu}=(\mu_1,\mu_2,\ldots, \mu_n),$ where
$\mu_1\le\mu_2\le\cdots\le \mu_n,$ not all equal, and ${\bf
\lambda}=(\lambda_1,\lambda_2,\ldots,\lambda_n),$ where
$\lambda_1\le\lambda_2\le\cdots\le \lambda_n,$ not all equal. Denote
\[\overline{\mu}=\frac1n\sum_{i=1}^n\mu_i,\ \ \overline{\lambda}=\frac1n\sum_{i=1}^n\lambda_i, \quad s^2(\mu)=\frac1n\sum_{i=1}^n(\mu_i-\overline{\mu})^2,\ \
s^2(\lambda)=\frac1n\sum_{i=1}^n(\lambda_i-\overline{\lambda})^2,
\]
and the order covariance \[{\rm
Cov}(\mu,\lambda)=\frac1n\sum_{i=1}^n\mu_i\lambda_i-\overline{\mu}\cdot\overline{\lambda}.\]
Then we have \begin{eqnarray*}\frac{1}{n-1}\le \frac{{\rm
Cov}(\mu,\lambda)}{s(\mu)\,s(\lambda)}. \end{eqnarray*} Equality
holds if and only if $($i\,$)$  $n=2$ or $($ii\,$)$ $\mu_1\le
\mu_2=\cdots=\mu_n$ and $\lambda_1=\cdots=\lambda_{n-1}\le
\lambda_n$ or $($iii\,$)$ $\mu_1= \cdots=\mu_{n-1}\le \mu_n$ and
$\lambda_1\le \lambda_2=\cdots=\lambda_n.$}
\medskip\\
\noindent{\bf Remark 2.} The main purpose of Benedetti's Inequality
here is to derive Corollary 2 above.\\
 \indent We define a compact
subset of ${\bf A}$ as follows:
\begin{eqnarray}{\bf A}_n=\Big\{(\lambda_1,\ldots,\lambda_n): \lambda_1\le \lambda_2\le\cdots\le \lambda_n, \
\small\sum_{i=1}^n\lambda_i=0,\ \sum_{i=1}^n\lambda_i^2=n-1\Big\}.
\end{eqnarray}

 \noindent{\bf Lemma 1.}  {\it Let $n\ge 3$ and let $Z_n$ defined in $(1)$ be
a feasible definite statistic on ${\bf A}$ with positive degree $1$
of homogeneity. Then there exist two positive constants $k<K$ such
that \begin{eqnarray}0<k\le\frac{Z_n}{S_n}\le K<\infty\ \  a.s.\
{\rm (almost\ surely)}.\end{eqnarray}} {\bf Proof.}  By the
homogeneity property, we rewrite
\begin{eqnarray*}Z_n&=&U(X_{(1)}-\overline{X}_n,X_{(2)}-\overline{X}_n,\ldots,X_{(n)}-\overline{X}_n)\\
&=&S_nU(\Lambda_1,\Lambda_2,\ldots,\Lambda_{n}),
\end{eqnarray*}
where
\begin{eqnarray*}\Lambda_i=\frac{X_{(i)}-\overline{X}_n}{S_n},\ \
i=1,2,\ldots,n,\end{eqnarray*} and
$(\Lambda_1,\Lambda_2,\ldots,\Lambda_{n})$ takes values
$(\lambda_1,\lambda_2,\ldots,\lambda_{n})$ in ${\bf A}_n$ (defined
in (13)) almost surely. Since the function $U$ is nonnegative,
definite and  continuous on ${\bf A}$ and the compact subset ${\bf
A}_n$ of ${\bf A}$ is bounded away from the origin, the range
$U({\bf A}_n)\subset(0,\infty)$ is a compact set as well. Hence,
there exist two positive constants $k<K$ such that $U({\bf
A}_n)\subset [k, K]\subset(0,\infty).$ This implies that
\[0<k\le U(\Lambda_1,\Lambda_2,\ldots,\Lambda_{n})\le K<\infty\ \  a.s.\]
Equivalently, (14) holds true. The proof is complete.\ \ ${\qed}$

As in the beginning of Section 2,  let $X_{(1)}\le X_{(2)}\le
\cdots\le X_{(n)}$ be the order statistics of a random sample $X_1,
X_2,\ldots, X_n$ of size $n\ge 3$ from a distribution $F$ which has
positive continuous density $f_X$ on $\mathbb R.$ The corresponding
realized order values are   $x_{(1)}\le x_{(2)}\le \cdots\le
x_{(n)}$ with the sample mean and the sample variance:
\[\overline{x}_n=\frac1n\sum_{i=1}^nx_{(i)},\quad
s_n^2=\frac{1}{n-1}\sum_{i=1}^n(x_{(i)}-\overline{x}_n)^2.\]

Following Hwang and Hu (1994b), we define the transformation
\begin{eqnarray*}T:(x_{(1)}, x_{(2)}, \ldots, x_{(n)})\longrightarrow
(t_1,t_2,\ldots,t_{n-2},w_1,w_2)\end{eqnarray*} by
\begin{eqnarray}
\begin{cases}\displaystyle
t_i=\Big[\frac{n-i+1}{(n-1)(n-i)}\Big]^{1/2}\Big[\frac{x_{(i)}-\overline{x}_n}{s_n}+\frac{1}{n-i+1}\sum_{k=1}^{i-1}\frac{x_{(k)}-\overline{x}_n}{s_n}\Big],\ \
1\le i\le n-2, \vspace{0.2cm}\\
w_1=\overline{x}_n, \vspace{0.1cm}\\
 w_2=s_n,\\
\end{cases}
\end{eqnarray}
where the summation is taken to be zero if $i=1.$ Then let
$f_{i}=1-\sum_{k=1}^it_k^2,$ where $1\le i\le n-2,$ and denote the
inverse transformation of $T$: \begin{eqnarray*}T^{-1}:
(t_1,t_2,\ldots,t_{n-2},w_1,w_2)\longrightarrow (x_{(1)}, x_{(2)},
\ldots, x_{(n)})\end{eqnarray*} through
\begin{eqnarray}
\begin{cases}\displaystyle
\frac{x_{(i)}-w_1}{w_2\sqrt{n-1}}=\Big[\frac{n-i}{n-i+1}\Big]^{1/2}\cdot t_i
-\sum_{k=1}^{i-1}\frac{t_k}{[(n-k)(n-k+1)]^{1/2}},\ \ 1\le i\le n-2, \vspace{0.2cm}\\\displaystyle
\frac{x_{(n-1)}-w_1}{w_2\sqrt{n-1}}=-\sum_{k=1}^{n-2}\frac{t_k}{[(n-k)(n-k+1)]^{1/2}}-[f_{n-2}/2]^{1/2}, \vspace{0.1cm}\\\displaystyle
\frac{x_{(n)}-w_1}{w_2\sqrt{n-1}}=-\sum_{k=1}^{n-2}\frac{t_k}{[(n-k)(n-k+1)]^{1/2}}+[f_{n-2}/2]^{1/2},\\
\end{cases}
\end{eqnarray}
where $\overline{x}_n=w_1$ and $s_n=w_2.$

 For $n\ge 3,$ define two more  subsets ${\bf D}_n$ and ${\bf R}_n$
of ${\mathbb R}^n$:
\begin{eqnarray}
{\bf D}_n&=&\{(x_{(1)}, x_{(2)},\ldots, x_{(n)}):\ x_{(1)}\le  x_{(2)}\le \cdots\le  x_{(n)}\},\\
{\bf R}_n&=&\Big\{(t_1,t_2,\ldots,t_{n-2},w_1,w_2):\ -1\le t_1\le
-1/(n-1),\nonumber\\
& &~~~\max\big\{\Big[\frac{n-k+2}{n-k}\Big]^{1/2}\cdot t_{k-1},\
-f_{k-1}^{1/2}\big\} \le t_k\le -f^{1/2}_{k-1}/(n-k),\nonumber\\
& &~~~\ 2\le k\le n-2,\ \ w_1\in{\mathbb R},\ \ w_2>0\Big\}.
\end{eqnarray}
Then we have the following:
\medskip\\
\noindent{\bf Lemma 2} (Hwang and Hu 1994b). {\it  For $n\ge 3,$ let
$T$ be the transformation  defined in $(15)$. Then $T$ establishes a
one-to-one correspondence between  the domain ${\bf D}_n$ and the
range ${\bf R}_n,$ defined in $(17)$ and $(18)$, respectively,
except for a set of $n$-dimensional Lebesgue measure zero.
Furthermore,  the absolute value of its Jacobian is
\[|J|=\sqrt{n}\,(n-1)^{(n-1)/2}\cdot w_2^{n-2}\cdot f_{n-2}^{-1/2}.
\]}
\noindent{\bf Lemma 3} (Hwang and Hu 1994b). {\it  Let $X\sim F$ be
a standard normal random variable. For $n\ge 3,$ define the
statistics
\begin{eqnarray}
T_i=\Big[\frac{n-i+1}{(n-1)(n-i)}\Big]^{1/2}\Big[\frac{X_{(i)}-\overline{X}_n}{S_n}+\frac{1}{n-i+1}\sum_{k=1}^{i-1}\frac{X_{(k)}-\overline{X}_n}{S_n}\Big],\ \
1\le i\le n-2.
\end{eqnarray}
Then $(T_1,T_2,\ldots,T_{n-2})$ has joint density
\begin{eqnarray}f(t_1,t_2,\ldots,t_{n-2})=\frac{n!\Gamma((n-1)/2)}{2\,\pi^{(n-1)/2}\cdot f_{n-2}^{1/2}},\quad
(t_1,t_2,\ldots, t_{n-2})\in {\bf B}_{n-2},
\end{eqnarray}
where $f_{n-2}=1-t_1^2-t_2^2-\cdots-t_{n-2}^2$ and ${\bf B}_{n-2}$
is a subset of ${\mathbb R}^{n-2}:$
\begin{eqnarray}
{\bf B}_{n-2}&=&\Big\{(t_1,t_2,\ldots,t_{n-2}):\ -1\le t_1\le
-1/(n-1),\nonumber\\
& &~~~\max\big\{\Big[\frac{n-k+2}{n-k}\Big]^{1/2}\cdot t_{k-1},\
-f_{k-1}^{1/2}\big\} \le t_k\le -f^{1/2}_{k-1}/(n-k),\nonumber\\
& &~~~\ 2\le k\le n-2\ \Big\}.
\end{eqnarray}}
\indent Comparing (15) and (19), we note that
$(T_1,T_2,\ldots,T_{n-2})$ has realized value
$(t_1,t_2,\ldots,t_{n-2})$ given in (15). Moreover, it follows from
the density (20) that
\[
\int_{{\bf B}_{n-2}}f_{n-2}^{-1/2}\prod_{i=1}^{n-2}{\rm d}t_i=\frac{2\,\pi^{(n-1)/2}}{n!\Gamma((n-1)/2)}<\infty,
\]
which will be used in the proof of the Theorem.
\medskip\\
 \noindent{\bf 4. Proofs of
main results}

\noindent {\bf Proof of the Theorem.} It suffices to prove the
necessity part. Namely, suppose that the sample mean
$\overline{X}_n$ is independent of the feasible definite statistic
\[Z_n=U(X_{(1)}-\overline{X}_n,X_{(2)}-\overline{X}_n,\ldots,X_{(n)}-\overline{X}_n)\]
on ${\bf A}$ with positive degree $p$ of homogeneity. Then we want
to prove that $F$ is normal. Without loss of generality, we may
assume that $p=1.$ Otherwise, we can consider instead  $Z_n^{1/p}$
which is also independent of $\overline{X}_n$ and is a feasible
definite statistic on ${\bf A}$ with positive degree $1$ of
homogeneity. Now, by the homogeneity property, rewrite, as in the
proof of Lemma 1,
\begin{eqnarray*}Z_n&=&U(X_{(1)}-\overline{X}_n,X_{(2)}-\overline{X}_n,\ldots,X_{(n)}-\overline{X}_n)\\
&=&S_nU(\Lambda_1,\Lambda_2,\ldots,\Lambda_{n}),
\end{eqnarray*}
where \begin{eqnarray}\Lambda_i=\frac{X_{(i)}-\overline{X}_n}{S_n},\
\ i=1,2,\ldots,n,\end{eqnarray} and
$(\Lambda_1,\Lambda_2,\ldots,\Lambda_{n})$ takes values
$(\lambda_1,\lambda_2,\ldots,\lambda_{n})$ in ${\bf A}_n$ (defined
in (13)) almost surely. It then follows from Lemma 1 that there
exist two positive constants $k<K$ such that
\begin{eqnarray}0<k\le \frac{Z_n}{S_n}=U(\Lambda_1,\Lambda_2,\ldots,\Lambda_{n})\le K<\infty\
\  a.s.\end{eqnarray}

Write the realized values of
\[
S_n=\frac{Z_n}{U(\Lambda_1,\Lambda_2,\ldots,\Lambda_{n})}\] as
\begin{eqnarray}s_n=\frac{z_n}{U(\lambda_1,\lambda_2,\ldots,\lambda_{n})}:=\frac{z_n}{U(\lambda)},
\end{eqnarray} where $\lambda=(\lambda_1,\ldots,\lambda_n).$

Recall that the order statistics $(X_{(1)}, X_{(2)}, \cdots,
X_{(n)})$ have joint density
\[
f(x_{(1)}, x_{(2)}, \ldots,
x_{(n)})=n!\prod_{i=1}^nf_X(x_{(i)}),\ \ x_{(1)}\le x_{(2)}\le \cdots\le
x_{(n)}.
\]
Next, consider the composition of two transformations $T$ (defined
in (15)) and $T^*$ defined by
\[T^*:(t_1,t_2,\ldots,t_{n-2},w_1,w_2)\longrightarrow (t_1,t_2,\ldots,t_{n-2},w_1,z_n),\]
where $w_1=\overline{x}_n,\  w_2=s_n,$  and $z_n$ is defined in
(24). The latter $T^*$  has Jacobian $J^*=1/U(\lambda)>0$ due to
(24). Hence, the composition $T^*\circ T$ has a Jacobian with
absolute value (by Lemma 2):
\begin{eqnarray*} |J|J^*&=&\sqrt{n}\,(n-1)^{(n-1)/2}\cdot
s_n^{n-2}\cdot f_{n-2}^{-1/2}\cdot\frac{1}{U(\lambda)}\\
&=&\sqrt{n}\,(n-1)^{(n-1)/2}\cdot
\Big(\frac{z_n}{U(\lambda)}\Big)^{n-2}\cdot f_{n-2}^{-1/2}\cdot\frac{1}{U(\lambda)}\\
&=&\sqrt{n}\,(n-1)^{(n-1)/2}\cdot
z_n^{n-2}\Big(\frac{1}{U(\lambda)}\Big)^{n-1}\cdot f_{n-2}^{-1/2},
\end{eqnarray*}
where $f_{n-2}=1-\sum_{i=1}^{n-2}t_i.$  Therefore, the statistics
$(T_1,T_2,\ldots,T_{n-2},\overline{X}_n,Z_n)$ have a joint density
on ${\bf R}_n$ (see (18)):
\[f(t_1,t_2,\ldots,t_{n-2},\overline{x}_n,z_n)=n!\prod_{i=1}^nf_X(x_{(i)})|J|J^*,
\]
where, by (16), (22)  and (24),
\[x_{(i)}=\overline{x}_n+\frac{z_n}{U(\lambda)}\cdot\lambda_i(t_1,t_2,\ldots,t_{n-2}),\quad i=1,2,\ldots,n.\]
This in turn implies that the joint density of
$(\overline{X}_n,Z_n)$ has the form (here, for simplicity, denote
$\overline{x}=\overline{x}_n$ and $z=z_n$):
\begin{eqnarray}
f(\overline{x},z)&=&n!\sqrt{n}\,(n-1)^{(n-1)/2}\cdot
z^{n-2}\int_{{\bf B}_{n-2}}\Big(\frac{1}{U(\lambda)}\Big)^{n-1}\nonumber\\
&  &\times\, f_{n-2}^{-1/2}\prod_{i=1}^n
f_X\Big(\overline{x}+\frac{z}{U(\lambda)}\cdot\lambda_i(t_1,t_2,\ldots,t_{n-2})\Big)\prod_{i=1}^{n-2}{\rm d}t_i,\ \ \overline{x}\in{\mathbb R},\ z\ge0,~~
\end{eqnarray}
where the set ${\bf B}_{n-2}\subset{\mathbb R}^{n-2}$ is defined in
(21).

Now we apply the independence condition on $\overline{X}_n$ and
$Z_n,$ and write the joint density in (25) as the products of the
densities of $\overline{X}_n$ and $Z_n:$
\begin{eqnarray}
f(\overline{x},z)=f_{\overline{X}_n}(\overline{x})\,f_{Z_n}(z),\quad \overline{x}\in {\mathbb R}\ {\rm and}\  z\ge 0.
\end{eqnarray}
Letting $\overline{x}=0$ in (25) and (26), we get the density of
$Z_n:$
\begin{eqnarray}
f_{Z_n}(z)&=&\frac{1}{f_{\overline{X}_n}(0)}\,f(0,z)\nonumber\\
&=& \frac{1}{f_{\overline{X}_n}(0)}\,n!\sqrt{n}\,(n-1)^{(n-1)/2}\cdot
z^{n-2}\int_{{\bf B}_{n-2}}\Big(\frac{1}{U(\lambda)}\Big)^{n-1}\nonumber\\
&  &\times\, f_{n-2}^{-1/2}\prod_{i=1}^n
f_X\Big(\frac{z}{U(\lambda)}\cdot\lambda_i(t_1,t_2,\ldots,t_{n-2})\Big)\prod_{i=1}^{n-2}{\rm d}t_i,\quad z\ge 0.
\end{eqnarray}
Plugging (25) and (27) in (26) and canceling the common term
$\sqrt{n}\,(n-1)^{(n-1)/2}\cdot z^{n-2}$, and then letting $z\to 0$
and  canceling the common integral term, we obtain the density of
$\overline{X}_n:$
\begin{eqnarray}
f_{\overline{X}_n}(\overline{x})=c\,(f_X(\overline{x}))^n,\quad \overline{x}\in {\mathbb R},
\end{eqnarray}
where $c=f_{\overline{X}_n}(0)/{f_X^n(0)}>0$ is a constant.

Combining (25) through (28), we finally have the integro-functional
equation\,:
\begin{eqnarray}
& &\int_{{\bf B}_{n-2}}\Big(\frac{1}{U(\lambda)}\Big)^{n-1}\cdot f_{n-2}^{-1/2}
\prod_{i=1}^n
f_X\Big(\overline{x}+\frac{z}{U(\lambda)}\cdot\lambda_i(t_1,t_2,\ldots,t_{n-2})\Big)\prod_{i=1}^{n-2}{\rm d}t_i\nonumber\\
&=& C\,(f_X(\overline{x}))^n\int_{{\bf B}_{n-2}}\Big(\frac{1}{U(\lambda)}\Big)^{n-1}\cdot f_{n-2}^{-1/2}
\prod_{i=1}^n
f_X\Big(\frac{z}{U(\lambda)}\cdot\lambda_i(t_1,t_2,\ldots,t_{n-2})\Big)\prod_{i=1}^{n-2}{\rm d}t_i,~~~
\end{eqnarray}
for all $\overline{x}\in{\mathbb R}$ and $z\ge 0,$ where
$C=1/f_X^n(0)>0$ is a constant.

It is seen that (29) is of the form of Anosov's integro-functional
equation (12), because
\[\int_{{\bf B}_{n-2}}\Big(\frac{1}{U(\lambda)}\Big)^{n-1}\cdot f_{n-2}^{-1/2}\,\prod_{i=1}^{n-2}{\rm d}t_i\ \ \in\ (0,\infty),
\]
\[\sum_{i=1}^n\sigma_i(t_1,t_2,\ldots,t_{n-2}):=\sum_{i=1}^{n}\frac{\lambda_i(t_1,t_2,\ldots,t_{n-2})}{U(\lambda)}=0,
\]
\[\sum_{i=1}^n\sigma_i^2(t_1,t_2,\ldots,t_{n-2})=\sum_{i=1}^{n}\Big[\frac{\lambda_i(t_1,t_2,\ldots,t_{n-2})}{U(\lambda)}\Big]^2=\frac{n-1}{U^2(\lambda)}\ \ \in\ (0,\infty).
\]
The last two conditions are required in (11).  We can check these
conditions by using  (22)--(24) and the remarks right after Lemma 3.
Besides, ${\bf B}_{n-2}, \overline{x}$ and $z$ here play the roles
of $\Phi, t$ and $s$ in (10) and (12), respectively. Therefore,
$f_X$ is normal by mimicking the proof of Anosov's theorem. This
completes the proof of the theorem.\ \ $\qed$
\medskip\\
\noindent {\bf Proof of Corollary 2.} Under the conditions on $a_i,\
i=1,2,\ldots,n,$  the base function
\[U(\lambda_1,\lambda_2,\ldots,\lambda_n)=\sum_{i=1}^na_i\lambda_i
\]
of (4) is  feasible definite on ${\bf A}$ with positive degree $1$
of homogeneity. To check $U\ge 0$ and the definiteness, write
$\overline{a}=\frac1n\sum_{i=1}^na_i$ and
$\overline{\lambda}=\frac1n\sum_{i=1}^n\lambda_i=0,$ then, by
Benedetti's inequality,
\[\sum_{i=1}^na_i\lambda_i=\sum_{i=1}^n(a_i-\overline{a})(\lambda_i-\overline{\lambda})\ge
\frac{1}{n-1}\Big[{\sum_{i=1}^n(a_i-\overline{a})^2}\Big]^{1/2}\Big[{\sum_{i=1}^n(\lambda_i-\overline{\lambda})^2}\Big]^{1/2}\ge 0.
\]
The LHS equals zero if and only if $\lambda_i=\overline{\lambda}=0$
for each $i,$ because $a_i$ are not all equal.\ \ $\qed$
\medskip\\
\noindent {\bf Proof of Corollary 3.} Under the conditions on $p,\
a_i,\ i=1,2,\ldots,n,$ it is easy to check that the base function
\[U(\lambda_1,\lambda_2,\ldots,\lambda_n)=\sum_{i=1}^na_i|\lambda_i|^p
\]
of (5) is  feasible definite on ${\bf A}$ with positive degree $p$
of homogeneity.\ \ $\qed$
\medskip\\
\noindent {\bf Proof of Corollary 4.} Under the conditions on $p,\
a_{ij},\ 1\le i,j \le n,$  it is easy to check that the base
function
\[U(\lambda_1,\lambda_2,\ldots,\lambda_n)=\sum_{i=1}^n\sum_{j=1}^na_{ij}|\lambda_i-\lambda_j|^p
\]
of (6) is  feasible definite on ${\bf A}$ with positive degree $p$
of homogeneity.\ \ $\qed$
\medskip\\
\noindent {\bf Proof of Corollary 6.} Under the conditions on
$a_{ij},\ 1\le i,j \le n,$  it is easy to check that the base
function
\[U(\lambda_1,\lambda_2,\ldots,\lambda_n)=\sum_{i=1}^n\sum_{j=1}^na_{ij}\lambda_i\lambda_j
\]
of (8) is  feasible definite on ${\bf A}$ with positive degree $2$
of homogeneity.\ \ $\qed$
\medskip\\
\noindent {\bf Proof of Corollary 7.} Under the conditions on $p,
q,\,a_{ij},\ 1\le i,j \le n,$  it is easy to check that the base
function
\[U(\lambda_1,\lambda_2,\ldots,\lambda_n)=\sum_{i=1}^n\sum_{j=1}^na_{ij}|\lambda_i|^p|\lambda_j|^q
\]
of (9) is  feasible definite on ${\bf A}$ with positive degree $p+q$
of homogeneity.\ \ $\qed$
\medskip\\
 \noindent{\bf 5. More results about the
sample range and Gini's mean difference}

For the case of sample size $n=2,$ a direct calculation shows the
interesting  relation between the sample variance and the sample
range: \[S_2^2=\frac12R_{2}^2.\] Hence if the sample mean
$\overline{X}_2$ and the sample range $R_{2}$ are independent, then
so are $\overline{X}_2$ and $S_2^2.$ This in turn implies that the
underlying distribution $F$ is normal. Therefore, in this case
($n=2$), we don't need to assume the smoothness conditions on the
distribution $F.$ Similarly, for Gini's mean difference in (7), we
have
\[G_2=R_2.
\]
Hence, if $\overline{X}_2$ and $G_2$ are independent, then $F$ is
normal.

In view of the proof of the Theorem, Lemma 1 plays a crucial role.
For the case of the sample range $R_n,$ we actually have,  by using
Benedetti's inequality, the following explicit bounds in the
inequality (14)\,:
\[\frac{\sqrt{2}}{\sqrt{n-1}}\le \frac{R_n}{S_n}\le \sqrt{2(n-1)}\ \ a.s.,\ \ n\ge 2.
\]
Further, a better lower bound ($\sqrt{2}$\,) can be obtained by
Lemma 4 and Proposition 1 below.
\medskip\\
\noindent{\bf Lemma 4.} {\it For any $n$ real numbers
$x_1,x_2,\ldots,x_n,$ we have the identity\,$:$
\[\sum_{i=1}^n(x_i-\overline{x}_n)^2=\frac{1}{2n}\sum_{i=1}^n\sum_{j=1}^n(x_i-x_j)^2=\frac1n\sum_{i<j}(x_i-x_j)^2,\ \ {\rm where }\ \ \overline{x}_n=\frac1n\sum_{i=1}^nx_i.
\]}
{\bf Proof.} Write $x_i-x_j=(x_i-\overline{x}_n)+
(\overline{x}_n-x_{j})$ and carry out the double summation.\ \
$\qed$
\medskip\\
 \noindent{\bf Proposition 1} (Range Inequality). {\it  Assume that  $x_{(1)}\le x_{(2)}\le \cdots\le x_{(n)}$ are not all equal, where $n\ge 2,$ and
denote
\[\overline{x}_n=\frac1n\sum_{i=1}^nx_{(i)}\ \ \ {\rm and}\ \ \ s^2_n=\frac{1}{n-1}\sum_{i=1}^n(x_{(i)}-\overline{x}_n)^2.\]
Then we have the range
inequality\,$:$\begin{eqnarray}{\sqrt{2}}\,s_n\le x_{(n)}-x_{(1)}\le
\sqrt{2(n-1)}\,s_n, \end{eqnarray} or, equivalently,
\begin{eqnarray}\sqrt{2}\le \frac{x_{(n)}-x_{(1)}}{s_n}\le \sqrt{2(n-1)}.
\end{eqnarray}}
{\bf Proof.} It follows from  Lemma 4 that
\begin{eqnarray}\sum_{i=1}^n(x_{(i)}-\overline{x}_n)^2=\frac1n\sum_{i<j}(x_{(i)}-x_{(j)})^2\le
\frac1n{n\choose
2}(x_{(n)}-x_{(1)})^2=\frac{n-1}{2}(x_{(n)}-x_{(1)})^2.
\end{eqnarray}
This proves the lower bound of $x_{(n)}-x_{(1)}$ in (30). To prove
the upper bound, we apply Cauchy--Schwarz inequality. More
precisely, take $y_{(1)}=-1, y_{(2)}=0, \ldots,
y_{(n-1)}=0,y_{(n)}=1.$ Then $\overline{y}_n=0$ and
\begin{eqnarray*}x_{(n)}-x_{(1)}&=&(x_{(n)}-\overline{x}_n)-(x_{(1)}-\overline{x}_n)=\sum_{i=1}^n(x_{(i)}-\overline{x}_n)(y_{(i)}-\overline{y}_n)\\&\le&
\Big({\sum_{i=1}^n(x_{(i)}-\overline{x}_n)^2}\Big)^{1/2}\Big({\sum_{i=1}^n(y_{(i)}-\overline{y}_n)^2}\Big)^{1/2}
=\sqrt{2(n-1)}\,s_n.
\end{eqnarray*}
This completes the proof.\ \ $\qed$

The result (31) can be used to derive the following inequality.
\medskip\\
\noindent{\bf Proposition 2.} {\it Under the conditions of Corollary
$4,$
\[(a_{1n}+a_{n1})
(\sqrt{2}\,)^p\le\frac{Z_n}{S_n^p}\le \Big(\max_{\substack{1\le i, \,j \le n\\
i\ne j}}{a_{ij}}\Big)\,n(n-1) (\sqrt{2(n-1)}\,)^p\ \ a.s.
\]In particular, for Gini's mean difference, we have
\begin{eqnarray}\frac{2\sqrt{2}}{n(n-1)}\le \frac{G_n}{S_n}\le \sqrt{2(n-1)}\ \
\ a.s.\end{eqnarray}}

\noindent{\bf Remark 3.} When $n=2,$  the upper and lower bounds in
(31) are  equal  and hence (31) becomes an identity (the same is
true for (33)). Moreover, for $n\ge 3,$ consider the symmetric case:
$-1,0,\ldots,0,1,$ then $\bar{x}_n=0,$ $x_{(1)}=-1,\ x_{(n)}=1,$
 and $s_n^2={2/(n-1)}.$ Hence, the right equality
in (31) also holds true in this case. This means that the upper
bound in (31) is sharp. However, for $n\ge 3,$ it follows from the
proof of Proposition 1 (see (32)) that $\sqrt{2}s_n=x_{(n)}-x_{(1)}$
if and only if all $x_{(i)}$ are equal. Therefore, under the
assumption that the  $x_{(i)}$ are not all equal, it is still
possible to improve the lower bound in (31) for the case $n\ge 3$
(see, e.g., Thomson 1955). On the other hand, when $n\ge 3,$ we have
$G_n<R_n$ {\it a.s.} by (7). So it is possible to improve both the
upper and lower bounds in (33) for the case $n\ge 3$ (see Barker
1983).
\medskip\\
\indent In view of the above observations and Corollaries 1 and 5,
we would like to pose the following conjectures. Equivalently, this
is to conjecture that  for the cases of the sample range and Gini's
mean difference, we don't need to assume the smoothness condition on
the underlying distribution.
\medskip\\
\noindent{\bf Conjecture 1.} {\it Let $X_1, X_2,\ldots, X_n$ be a
random sample of size $n\ge 3$ from {\it any} distribution $F$ on
$\mathbb R.$  If the sample mean $\overline{X}_n$ and the sample
range $R_{n}$ in $(3)$ are independent, then $F$ is a normal
distribution $($including the degenerate case$).$}
\medskip\\
\noindent{\bf Conjecture 2.} {\it Let $X_1, X_2,\ldots, X_n$ be a
random sample of size $n\ge 3$ from {\it any} distribution $F$ on
$\mathbb R.$  If the sample mean $\overline{X}_n$ and Gini's mean
difference $G_n$ in $(7)$  are independent, then $F$ is a normal
distribution $($including the degenerate case$).$}
\medskip\\
\noindent{\bf 6. Discussions}

Let $X_1, X_2,\ldots, X_n$ be a random sample of size $n\ge 2$ from
normal distribution $N(\mu,\sigma^2).$ As before, denote the sample
mean and the sample variance by $\overline{X}_n$ and $S_n^2,$
respectively. Then, it is known that the sampling distribution of
the useful statistic
\[T_{n-1}=\frac{\overline{X}_n-\mu}{S_n/\sqrt{n}}
\]
is  Student's $t$-distribution with $n-1$ degrees of freedom. The
derivation of the sampling distribution heavily depends on the
independence of $\overline{X}_n$ and $S_n^2$ (or ${S_n}).$
Analogously, instead of ${S_n/\sqrt{n}},$ let $Z_n$ be another
feasible definite statistic with positive degree 1 of homogeneity
(e.g., the sample range or Gini's mean difference; both are linear
functions of order statistics), then it is independent of
$\overline{X}_n.$ If one can carry out the sampling distribution of
the new statistic
\[
T^*_{n-1}=\frac{\overline{X}_n-\mu}{Z_n},
\]
 it will be useful in statistical inference, such as finding the
confidence interval (CI) of the mean $\mu$ when $\sigma$ is unknown.
This allows us to compare the CIs obtained through $T_{n-1}$ and
$T^*_{n-1}.$ Besides, the characterization results  provided here
might be  useful in hypothesis testing, namely, testing the
normality assumption via the independence of the sample mean
 and some feasible definite statistics.
 \medskip\\
 \noindent{\bf
Acknowledgments.} The authors would like to thank the Chief Editor,
Associate Editor and two Referees for helpful comments and
suggestions, which improve the presentation of the manuscript. We
also thank Professor Jordan Stoyanov for the complete information of
the references Anosov (1964) and Laha (1956).
\medskip\\
\noindent{\bf References}
\begin{description}
\item Anosov, D.\,V. (1964).
On an integral equation from statistics. {\it Vestnik of
Leningrad State University, Series III Mathematics Mechanics
Astronomy}, {\bf 19}, 151--154. (In Russian)

\item Barker, L. (1983). On Gini's mean difference and the sample standard deviation.
{\it Communications in Statistics -- Simulation and
Computation}, {\bf 12}, 503--505.

\item Benedetti, C. (1957).  Di alcune disuguaglianze collegate al campo di variazione di
indici statistical coefficients. {\it Metron}, {\bf  18},
102--125.

\item Benedetti, C. (1995).
On some inequalities related to the range of some statistical
indices. {\it Metron}, {\bf 52}, 189--215.

\item David, H.\,A. and Nagaraja, H.\,N. (2003). {\it Order Statistics}, 3rd edn. New Jersey: Wiley.

\item Geary, R.\,C. (1936). The distribution of ``Student's" ratio for non-normal samples.
{\it Supplement to the Journal of the Royal Statistical
Society}, {\bf 3}, 178--184.

\item Georgescu-Roegen, N. (1959).
On the extrema of some statistical coefficients. {\it Metron},
{\bf 19}, 38--45.

\item Hwang, T.-Y. and Hu, C.-Y. (1994a).
The best lower bound of sample correlation coefficient with
ordered restriction. {\it Statistics $\&$ Probability Letters},
{\bf 19}, 195--198.

\item Hwang, T.-Y. and Hu, C.-Y. (1994b).  On the joint
distribution of Studentized order statistics.  {\it Annals of
the Institute of Statistical Mathematics},  {\bf 46}, 165--177.

\item Hwang, T.-Y. and Hu, C.-Y. (2000). On some characterizations of population distributions. {\it  Taiwanese Journal of Mathematics}, {\bf 4}, 427--437.

\item Kagan, A.\,M,  Linnik, Yu.\,V.  and Rao, C.\,R. (1973). {\it Characterization Problems in
 Mathematical Statistics.}  New York: Wiley.

 \item Kawata, T. and Sakamoto, H. (1949). On the characterization of the normal population by
 the independence of the sample mean and the sample variance.
 {\it Journal of the Mathematical Society of Japan}, {\bf 1},
 111--115.

 \item Laha, R.\,G. (1956). On stochastic independence of a homogeneous quadratic statistic and of the mean. {\it Vestnik
of Leningrad State University, Series III Mathematics Mechanics
Astronomy}, {\bf 11}, 25--32. (In Russian)

\item Lukacs, E. (1942). A characterization of the normal distribution. {\it The Annals of Mathematical Statistics}, {\bf 13}, 91--93.

 \item Rohatgi, V.\,K. (1976). {\it An Introduction to Probability Theory and Mathematical Statistics.} New York: Wiley.

 \item Sarria, H. and Mart{i}nez, J.\,C. (2016). A new proof of the Benedetti's inequality and some applications to perturbation to real
 eigenvalues and singular values. {\it Boletin de
 Matem\'aticas}, {\bf 23}, 105--114.

\item Thomson, G.\,W. (1955). Bounds for the ratio of range to standard deviation. {\it Biometrika}, {\bf 42}, 268--269.

 \item Zinger, A.\,A. (1951). On independence samples from normal populations. {\it Uspekhi Matematicheskikh Nauk}, {\bf 6}, 172--175.
(In Russian)

\end{description}
\end{document}